\documentclass[11pt]{article}

\usepackage{latexsym,amsmath,amssymb,mathrsfs,graphicx,hyperref,mathpazo,soul,color,amsthm}
\hypersetup{pdftex,colorlinks=true,linkcolor=blue,citecolor=blue,filecolor=blue,urlcolor=blue}

\usepackage{url}
\usepackage{hyperref}

\renewcommand{\le}{\leqslant}
\def\tr{{\raise0pt\hbox{$\scriptscriptstyle\top$}}}

\newtheorem{theorem}{Theorem}[section]

\newtheorem{proposition}[theorem]{Proposition}

\newtheorem{problem}[theorem]{Problem}
\newtheorem{definition}[theorem]{Definition}
\newtheorem{question}[theorem]{Open Question}

\title{\vspace{-1.5 cm}{\bf On the smallest open Diophantine equations}}

\author{Bogdan Grechuk\footnote{School of Computing and Mathematical Sciences, University of Leicester, LE1 7RH, UK; bg83@leicester.ac.uk}}

\begin{document}

\maketitle

\begin{abstract}
This paper reports on the current status of the project in which we order all polynomial Diophantine equations by an appropriate version of ``size'', and then solve the equations in that order. We list the ``smallest'' equations that are currently open, both unrestricted and in various families, like the smallest open symmetric, 2-variable or 3-monomial equations. All the equations we discuss are amazingly simple to write down but some of  them seem to be very difficult to solve. 
\end{abstract}


\section{Introduction.}

In many areas of mathematics researchers order instances of difficult problems by some natural parameters and try to solve the problems at least for ``small'' instances, while in other areas such systematic approach is currently missing. To illustrate this point, let us look at some recent breakthroughs in knot theory and in number theory. 
%
%In math we often order instances of a problem by some natural parameter.  We will later apply
%this approach to diophantine equations; however, for now we will give two examples where this
%was done with great success.

One of the most amazing recent results in knot theory is the 2020 theorem of Lisa Piccirillo, stating that the Conway knot is not slice \cite{piccirillo2020conway}. Here, a knot is called (smoothly) \emph{slice} if it bounds a smoothly embedded $2$-dimensional disk in $4$-dimensional space, and the Conway knot is a certain knot labeled 11n34 in the standard knot tables. Here, ``$11$'' refers to the fact the the diagram of this knot has $11$ crossings. The question whether the Conway knot is slice or not has been opened for about $50$ years, and has attracted a lot of attention because \emph{all other knots of under $13$ crossings were classified according to whether they are slice knots or not, and the Conway knot has been the only exception}.

One of the most amazing results in number theory is the Andrew Wiles' proof \cite{wiles1995modular} of Fermat Last Theorem. The Theorem states that exponential Diophantine equation $x^n+y^n=z^n$ has no non-trivial integer solutions for $n\geq 3$, has been an open question for many centuries, and has attracted a lot of attention because... Fermat claimed that he has a beautiful proof of this fact which he never wrote down, and many generations of mathematicians were trying to find the missing proof.

A more recent example of an amazing result in the area of Diophantine equations is the 2021 paper of Booker and Sutherland \cite{booker2021question}, which studies equations in the form 
\begin{equation}\label{eq:BookSuth}
x^3+y^3+z^3=k.
\end{equation} 
Note that \eqref{eq:BookSuth} has no solution if $k\equiv \pm 4 \pmod 9$. It is known that there are infinitely many integer solutions for $k=1$ and $2$. For $k=3$, \eqref{eq:BookSuth} has the obvious solutions $(1,1,1)$ and permutations of $(4,4,-5)$. Miller and Woolett \cite{miller1955solutions} asked if there is always a solution when $k\not\equiv \pm 4 \pmod 9$. Mordell \cite{mordell1953integer} asked if there are other solutions when $k=3$. Brooker and Sutherland showed that
\begin{itemize}
\item[(i)] $x^3+y^3 + z^3 = 42$ has a solution (combined with prior work this shows
that for $|k|\le 100$, $k\not\equiv \pm 4 \pmod 9$, there is always a solution), and
\item[(ii)] there is another solution to $x^3+y^3+z^3=3$, namely:
$
x=569,936,821,221,962,380,720, \, y=-569,936,821,113,563,493,509, \, z=-472,715,493,453,327,032.
$ 
\end{itemize}
%This equation has attracted a lot of attention because... Mordell \cite{mordell1953integer} asked in 1953 whether it has any integer solutions other than $(1,1,1)$ and permutations of $(4,4,-5)$.

As you can see, in knot theory researchers ordered all knots by a natural parameter, number of crossings, and then try to answer questions of interest \emph{systematically} for all knots with the given number of crossings. A similar approach is taken in many other areas of mathematics: if a problem is difficult or undecidable in general, researchers order the instances of the problem in some natural way, and then try to ``solve'' at least ``small'' instances. For example, the Halting problem has been solved for all $2$-symbol Turing machines with $2$, $3$ and $4$ states, so the next open case is $5$-state machines. In contrast, Diophantine equations are not studied in order, and often attract interest just because some time ago a famous mathematician wrote down an equation and asked to solve it. 
%Another reason for a Diophantine equation to attract attention of number theorists is when the equation arise naturally in some application   

One reason for this is that it is not obvious how to order all Diophantine equations is a natural way. It is easy to arrange in order some specific infinite families of equations. For example, Fermat Last Theorem can be treated as one exponential equation or as an infinite family of polynomial Diophantine equations ordered by the value of parameter $n$. Similarly, equations \eqref{eq:BookSuth} can be naturally ordered by parameter $k$. However, it is less easy to order all polynomial Diophantine equations, that is, all equations of the form
\begin{equation}\label{eq:diofgen}
P(x_1,\dots, x_n)=0,
\end{equation}
where $P$ is a polynomial with integer coefficients. For this, we need to assign to every Diophantine equation a ``size'' parameter, such that for any bound $B$ there is only a finite number of equations of size at most $B$. Many natural notions of ``size'' or ``height'' do not satisfy this condition. For example, if we define the height of the equation \eqref{eq:diofgen} as the maximum (or sum) of absolute values of coefficients of $P$, then there are infinitely many equations of height $1$, e.g. $x^n=0$, $n=1,2,\dots$.

It is natural to define the ``size'' of equation \eqref{eq:diofgen}   
as a sum of sizes of monomials of $P$, so it is left to define the size of a monomial $ax_1^{k_1}\dots x_n^{k_n}$. If we consider an equation as an input to a computer program which then tries to solve it, then standard way to measure the size of the input is the number of bits needed to describe it. For simplicity, let us assume that we do not use the power symbol, and write $x_1^{k_1}$ as $x_1x_1\dots x_1$ ($k_1$ times) and so on. Then we need $d=k_1+\dots+k_n$ symbols to write $x_1^{k_1}\dots x_n^{k_n}$, where $d$ is the degree of the monomial. We also need about $\log_2 |a|$ symbols to write the coefficient $a$ in binary, so we can define the length of the monomial as $l=\log_2 |a|+d$. This is not an integer, but ordering the monomials by $l$ is equivalent to ordering them by $H=2^l=|a|2^d$, which is an integer. Hence, let us define the size of equation \eqref{eq:diofgen} as
\begin{equation}\label{eq:Hdef}
H(P)=\sum_{i=1}^k |a_i| 2^{d_i},
\end{equation}
provided that polynomial $P$ consists of $k$ monomials with integer coefficients $a_1, \dots, a_k$ and degrees $d_1, \dots, d_k$, respectively. For example, for the equation \eqref{eq:BookSuth} with $k=3$ we have $H=2^3+2^3+2^3+3=27$. This notion of size has been suggested by anonymous mathoverflow user Zidane who asked in 2018 what is the smallest open Diophantine equation, see \cite{Z2018}. Note that $H$ is always a non-negative integer, and, for any $B$, there is a finite number of equations with $H \leq B$. Hence, we may list all equations with sizes $H=0,1,2,\dots$, try to solve them in this order, and report the smallest equations we cannot solve. This is exactly the purpose of this paper.

In \cite{grechuk2021diophantine}, this project is implemented for the Hilbert 10th problem, that is, the problem of determining whether an equation has any integer solution. In this paper, we consider more general problem of determining all integer solutions, but also re-iterate the open questions posted in \cite{grechuk2021diophantine}. Section \ref{sec:polfam} considers the most general problem of solving the equations completely, including the description of the solution set if it is infinite. Section \ref{sec:2varall} considers this problem for the $2$-variable equations only. In Section \ref{sec:fin}, we consider an easier problem of determining whether the solution set is finite, and if so, list all the solutions. In Section \ref{sec:Hilbert}, we study an even easier problem of determining whether a given equation has any integer solution or not. For each of these problems, we have identified the smallest equations for which the problem is non-trivial, and challenge the readers to try to solve these equations. 
%, but the intersection of the material is still substantial. One may consider this paper as a short version of \cite{grechuk2021diophantine}, written for the open problem column. In particular, we refer to \cite{grechuk2021diophantine} for all the proofs.
 
\section{Describing all solutions: polynomial families}\label{sec:polfam}

We have exactly one equation $0=0$ of size $H=0$, and exactly two equations $\pm 1=0$ of size $H=1$. More generally, for every $H>0$, we have two equations $\pm H=0$ with no variables and no solutions. If we do not count these as ``equations'' and insist that every equation should have at least one variable, then the smallest equations are $\pm x=0$ of size $H=2$ and the nest smallest are $\pm x \pm 1=0$ of size $H=3$. All these equations are trivial to solve. More generally, for any equation in one variable
\begin{equation}\label{eq:onevar}
a_m x^m + \dots + a_1 x + a_0 = 0
\end{equation}
with integer coefficients $a_m,\dots, a_0$, any non-zero integer solution must be a divisor of $a_k$, where $k$ is the smallest integer such that $a_k\neq 0$. This gives an algorithm to list all integer solutions of \eqref{eq:onevar}. Moreover, this can be done in polynomial time, see \cite{cucker1999polynomial}. 

Further, we will consider only equations in at least two variables. The smallest such equations are $\pm x \pm y = 0$ and $\pm xy=0$ of size $H=4$. To avoid solving essentially the same equations multiple times, we call two equations equivalent if one can be transformed into another after multiplication by $-1$ and/or substitutions in the form $x_i \to -x_i$. Obviously, it suffices to consider only one equation from each equivalence class. With this convention, the only equations of size $H=4$ we need to consider are 
\begin{equation}\label{eq:xpy}
x + y = 0
\end{equation}
and 
\begin{equation}\label{eq:xy}
xy = 0.
\end{equation}
Both these equations have infinitely many integer solutions, which can be presented in a parametric form. For equation \eqref{eq:xpy}, the solutions are $(x,y)=(u,-u)$, where $u$ is an integer parameter. This is a special case of a polynomial family, as defined below.

\begin{definition}
We say that a subset $S \subset {\mathbb Z}^n$ is a polynomial family if there exist polynomials $P_1, \dots, P_n$ is $k$ variables $u_1,\dots, u_k$ such that $(x_1,\dots, x_n)\in S$ if and only if there exists integers $u_1,\dots, u_k$ such that $x_i=P_i(u_1,\dots, u_k)$, $i=1,\dots, n$. 
\end{definition}

In this terminology, the solution set of the equation \eqref{eq:xpy} is a polynomial family with $k=1$ parameter $u$, $P_1(u)=u$ and $P_2(u)=-u$. More generally, the solution set of any equation in the form
\begin{equation}\label{eq:xnpQ}
x_n + Q(x_1,\dots,x_{n-1}) = 0
\end{equation}
can be represented as polynomial family $x_i=u_i$, $i=1,\dots, n-1$, $x_n=-Q(u_1,\dots,u_{n-1})$.

The solution set to the equation \eqref{eq:xy} is $(x,y)=(0,u)$ or $(x,y)=(u,0)$ for any integer $u$. Note that it is not a polynomial family but a union of two polynomial families. Obviously, if we can represent the solution set as a union of any finite number of polynomial families, then we can classify the equation as ``solved''. Note that if the solutions sets of equations $P_1=0$ and $P_2=0$ are finite unions of polynomial families, then the same is true for the equation
\begin{equation}\label{eq:P1P2}
P_1 \cdot P_2 = 0.
\end{equation}  
From now on, we will exclude the equations of the forms \eqref{eq:xnpQ} and \eqref{eq:P1P2} from further analysis.

For $H=5$, we would like to mention equation
\begin{equation}\label{eq:2xp1}
2x +1 = 0.
\end{equation}
It is a one-variable equation \eqref{eq:onevar}, and we have already discussed and excluded all such equations. However, it is interesting as the smallest equation for which the left-hand side is always odd and therefore cannot be equal to $0$. More generally, if there exists an integer $m\geq 2$ such that $P(x_1,\dots,x_n)$ is never divisible by $m$, then equation \eqref{eq:diofgen} has no integer solutions. Hence, we may exclude such equations.

The only non-excluded equation of size $H=5$ is  
\begin{equation}\label{eq:xyp1}
xy +1 = 0
\end{equation}
whose solutions are $(x,y)=(1,-1)$ and $(x,y)=(-1,1)$. There are no other integer  solutions because in every potential solution $(x,y)$ both $x$ and $y$ must be a divisor of $1$. More generally, any equation of the form
\begin{equation}\label{eq:P1P2pc}
P_1 \cdot P_2 = c
\end{equation} 
for some integer $c$ can be solved by enumerating divisors $d$ of $c$, and, for every divisor, solve the system of equations $P_1=d, P_2=c/d$. From now on, we will exclude the equations of the form \eqref{eq:P1P2pc} solvable by this method. Families \eqref{eq:onevar}, \eqref{eq:xnpQ} and \eqref {eq:P1P2pc} cover all the equations of size $H \leq 7$. 

For $H=8$, there are some equations not considered so far, for example, equation
\begin{equation}\label{eq:x2py^2}
x^2 + y^2 = 0,
\end{equation}
whose only integer solution is $x=y=0$ because this is the only real solution. More generally, we can exclude all equations \eqref{eq:diofgen} in $n$ variables whose set of real solutions is a bounded region in ${\mathbb R}^n$, because any bounded region has at most a finite number of integer points, and we may use the direct substitution to check which of these points are the solutions to \eqref{eq:diofgen}.

Another equation not excluded so far is
\begin{equation}\label{eq:xyp2z}
xy+2z=0.
\end{equation}
With $z'=2z$, \eqref{eq:xyp2z} reduces to the equation $xy+z'=0$ with extra condition that $z'$ is even. The solution to the latter is a polynomial family $(x,y,z')=(u,v,-uv)$. By considering all possible parities of $u,v$, we can represent it as the union of $4$ polynomial families (i) $(x,y,z')=(2u',2v',-4u'v')$, (ii) $(x,y,z')=(2u',2v'+1,-2u'(2v'+1))$, (iii) $(x,y,z')=(2u'+1,2v',-(2u'+1)2v')$ and (iv) $(x,y,z')=(2u'+1,2v'+1,-(2u'+1)(2v'+1))$. In each case, the parity of $z'$ is determined: $z'$ is even in cases (i)-(iii) but odd in case (iv). By substituting $z=z'/2$ in formulas (i)-(iii), we obtain parametric solutions to \eqref{eq:xyp2z}.

Exactly the same argument implies, more generally, the following observation.  

\begin{proposition}\label{prop:linsub}
If the solution set to the equation $P(x_1,\dots,x_n)=0$ is a finite union of polynomial families, then so is the solution set to the equation 
$$
P(a_1x_1+b_1, \dots, a_nx_n+b_n)=0,
$$
where $a_1,\dots,a_n$ are non-zero integers and $b_1,\dots,b_n$ are arbitrary integers.
\end{proposition}

Proposition \ref{prop:linsub} allows us to do arbitrary linear substitutions in the form $ax_i+b_i \to x'_i$. For example, any equation of the form
%Because $xy$ is even, then either $x$ or $y$ must be even. In the first case, let us write $x=2u$ for some integer $u$, while in the second case $y=2v$ for some integer $v$. Hence, we get two families of integer solutions $(x,y,z)=(2u,v,-uv)$ and $(x,y,z)=(u,2v,-uv)$ for some integers $u,v$.
%More generally, any equation of the form
\begin{equation}\label{eq:axnpQ}
a x_n + Q(x_1,\dots,x_{n-1}) = 0,
\end{equation}
where $a\neq 0$ is an integer, 
%can be solved by enumerating all possible $|a|^{n-1}$ remainders $x_1, \dots, x_{n-1}$ can give after division by $a$, and in each case representing $x_i=au_i+r_i$, $i=1,\dots, n-1$, with $0\leq r_i < |a|$. Then each case when $Q(au_1+r_1, \dots, au_{n-1}+r_{n-1})$ is divisible by $a$ lead to a polynomial family of the solutions of \eqref{eq:axnpQ}. 
reduces to \eqref{eq:xnpQ} by substitution $a x_n \to x'_n$ and can therefore be excluded from further analysis.

All equations we have considered so far are trivial. The smallest equation that deserves to be given to students as an exercise is the equation
\begin{equation}\label{eq:x2pyz}
x^2-yz=0.
\end{equation} 
The question is, of course, how to represent all solutions as a polynomial family. To solve this exercise, the students should note that any integers $y$ and $z$ can be represented as $y=uv^2$ and $z=u'w^2$ with $u, u'$ square-free. Now, for $yz$ to be a perfect square we must have $u=u'$, hence $z=uw^2$, and then $x=uvw$. Conversely, for any $u,v,w$ (not necessarily square-free), the triple $(x,y,z)=(uvw,uv^2,uw^2)$ is a solution to \eqref{eq:x2pyz}.

Another easy exercise is to write as a polynomial family the set of all solutions to the equation
\begin{equation}\label{eq:xypzt}
xy-zt=0.
\end{equation} 
To solve it, denote $u_1$ the greatest common divisor of $x$ and $z$ (assuming that $x$ and $z$ are not both $0$). Then $x=u_1u_2$ and $z=u_1u_3$, where $u_2$ and $u_3$ are co-prime integers. Then $u_2y=u_3z$ implies that $z$ is divisible by $u_2$ and can be written as $z=u_2u_4$ for integer $u_4$. Then $y=u_3u_4$, and we obtain a polynomial family of solutions $(x,y,z,t)=(u_1u_2,u_3u_4,u_1u_3,u_2u_4)$. It is left to note that the remaining case $x=z=0$ is also covered by this parameterization.
%The answer is $(x,y,z,t)=(u_5u_1u_2, u_5u_3u_4, u_5u_1u_3, u_5u_2u_4)$ for integers $u_1,u_2,u_3,u_4,u_5$. It is trivial to check that the given family satisfies \eqref{eq:xypzt}. We leave it to the reader to show that, conversely, all solutions to \eqref{eq:xypzt} are covered by this parameterization.
This finishes the analysis of all equations of size $H \leq 8$.

With $H=9$, the problem suddenly jumps from student-exercise level to research level. The question whether the set of integer solutions of the equation
\begin{equation}\label{eq:xypztp1}
xy-zt=1
\end{equation}  
is a polynomial family has been first asked by Skolem \cite{skolem1938diophantische} in the 1930's, remained open for over $70$ years, and has been answered by Vaserstein \cite{vaserstein2010polynomial} in 2010, who proved that it is indeed a polynomial family with $46$ parameters. As a corollary of this result, Vaserstein also showed that, for any integer $c$, the solution set of the equation 
\begin{equation}\label{eq:x2pyzpc}
x^2-yz=c
\end{equation} 
is the union of a finite number of polynomial families. In particular, this covers the equations $x^2-yz=\pm 1$, and finishes the analysis of all equations of size $H \leq 9$. 

Vaserstein also proved that, for any $c$, the solution set of the equation 
\begin{equation}\label{eq:xypztpc}
x_1x_2+x_3x_4+Q(x_5,\dots,x_n)=c,
\end{equation} 
where $Q$ is a quadratic form, is the union of a finite number of polynomial families. For example, with $Q=0$ and $c=2$, this covers the equation $xy+zt=2$ of size $H=10$. Vaserstein's solution to \eqref{eq:xypztp1}-\eqref{eq:xypztpc} immediately allows to solve many other equations. For example, let $x,y,z$ be a solution to the equation
\begin{equation}\label{eq:x2pxpyz}
x^2+x-yz=0
\end{equation} 
of size $H=10$. 
Because $x^2+x=x(x+1)$, the prime factors of $y$ are distributed between $x$ and $x+1$, hence $y$ can be written as $y=ab$, where $a$ and $b$ are divisors of $x$ and $x+1$, respectively. With $c=\frac{x}{a}$ and $d=\frac{x+1}{b}$, we get $db-ac=(x+1)-x=1$, which is exactly \eqref{eq:xypztp1}, hence the set of such $(a,b,c,d)$ is a polynomial family. Then $(x,y,z)=(ac,ab,cd)$ is also a polynomial family. 

Alternatively, we may rewrite \eqref{eq:x2pxpyz} as $(2x+1)^2-1-4yz=0$. This equation reduces to $(x')^2-y'z=1$ after substitutions $x'=2x+1$ and $y'=4y$, which we can do by Proposition \ref{prop:linsub}. But the last equation is of the form \eqref{eq:x2pyzpc} and has been solved in \cite{vaserstein2010polynomial}.

More generally, any equation of the form
\begin{equation}\label{eq:ax2pbxpcpdyz}
ax^2+bx+c-dyz=0, \quad\quad a\neq 0
\end{equation}
can be rewritten as $4adyz=4a^2x^2+4abx+4ac = (2ax+b)^2 + 4ac-b^2$, which is equation of the form \eqref{eq:x2pyzpc} with $2ax + b \to x'$ and $4ady \to y'$. This allows us to conclude that the solution set to \eqref{eq:ax2pbxpcpdyz} is a finite union of polynomial families by Proposition \ref{prop:linsub} and exclude such equations from further analysis. More generally, we may exclude all equations that, after linear substitutions as in Proposition \ref{prop:linsub}, can be reduced to the already solved ones.
%equation with smaller $H$, or from already excluded family.

After this, the only remaining equations of size $H\leq 11$ are trivial two-variable quadratic equations like $y^2=x^2+x$ or $y^2=x^2+x\pm 1$. If $|x|>2$, then 
$
(|x|-1)^2 < x^2+x\pm 1 < (|x|+1)^2,
$
and also $x^2+x\pm 1 \neq |x|^2$, hence $x^2+x\pm 1$ cannot be a perfect square. By checking cases $|x|\leq 2$, we may easily list all the solutions. 

For $H=12$, we meet similar equations $y^2=x^2+x\pm 2$ and $y^2=x^2+2x$, that can be solved by exactly the same method. 
%A little bit different is the equation
%$$
%y^2+y=x^2+x,
%$$  
%which, after multiplication by $4$ and adding $1$, can be reduced to $(2y+1)^2=(2y+1)^2$, from which we derive two families of solutions $(x,y)=(u,u)$ and $(x,y)=(u,-u-1)$. 
A (trivial but) notable equation is
\begin{equation}\label{eq:2x2my2}
y^2 = 2x^2
\end{equation}
whose only integer solution is $(x,y)=(0,0)$. If it would have a non-zero integer solution $(x,y)$, then rational number $t=\frac{y}{x}$ would satisfy $t^2=2$, hence the non-existence of non-zero integer solutions to \eqref{eq:2x2my2} is equivalent to the irrationality of $\sqrt{2}$, a famous old problem with rich history. Similarly, equation
$
x^2+xy-y^2=0
$
reduces to non-existence of rational solutions to $t^2+t-1=0$. Next, equation
$$
2xy+x+y=0
$$
reduces to the question for which $x$ the ratio $\frac{x}{2x+1}$ can be an integer. It is easy to see that this is possible only for $x=-1$ and $x=0$, leading to the solutions $(x,y)=(-1,-1)$ and $(x,y)=(0,0)$. The more general equation
\begin{equation}\label{eq:liny}
Q(x)y+R(x)=0
\end{equation} 
reduces to the question when $\frac{R(x)}{Q(x)}$ is an integer, which can be easily answered in full generality, see \cite{grechuk2021diophantine}. Finally, all integer solutions to the equation $x^3-y^2=0$ are given by $(x,y)=(u^2,u^3)$.

This finishes the analysis of $2$-variable equations of size $H=12$. Among the $3$-variable ones, the most famous is the equation
\begin{equation}\label{eq:Pyth}
x^2+y^2=z^2
\end{equation}
whose integer solutions are known as Pythagorean triples. A standard approach to solving this equation is noting that if $z=0$ then $x=y=0$ and otherwise the equation can be written as $(x/z)^2+(y/z)^2=1$ and reduces to finding rational points on the circle $x^2+y^2=1$. To find them, we can choose any one rational point, say $(x,y)=(1,0)$, and draw all possible lines through this point with rational slope $k$. Any such line intersects the circle at $(1,0)$ and in another point which can be easily seen to be a rational point. This way we get parameterization of all rational points with rational parameter $k$. From this, we can easily write down all integer solutions to \eqref{eq:Pyth}. The answer is the union of two polynomial families $(x,y,z)=(w(u^2-v^2),2wuv,w(u^2+v^2))$ and $(x,y,z)=(2wuv,w(u^2-v^2),w(u^2+v^2))$. 

We continue our analysis of equations of size $H=12$ by discussing equations with two monomials, which are quite easy. For example, if $(x,y,z)$ satisfy the equation $x^2 y-z^2 = 0$, then $y$ must be a perfect square, hence the solution is $(x,y,z)=(u,v^2,uv)$. As another example, consider equation
$$
x^3-y z = 0.
$$
We can write $y=u_1^3u_2^2u_3$ and $z=u_4^3u_5^2u_6$ with $u_2,u_3,u_5,u_6$ co-prime and square-free. Then $yz$ is a perfect cube only if $u_5=u_3$ and $u_6=u_2$, resulting in the answer $(x,y,z)=(u_1u_2u_3u_4,u_1^3u_2^2u_3, u_4^3u_3^2u_2)$. A bit more difficult is the equation
$$
x^2 y-z t = 0
$$
The cases $y=0$ and $z=0$ can be considered separately, so we assume that $yz \neq 0$. Let $u_1=\text{gcd}(y,z)$, so that $y=u_1u_2$ and $z=u_1v_2$ with $u_2,v_2$ coprime. Then $x^2u_2 = t v_2$. Hence $t$ is divisible by $u_2$, let us write $t=u_2v_1$. Then $x^2=v_1v_2$, but this is the equation \eqref{eq:x2pyz}. Its solution is $x=u_3u_4u_5$, $v_1=u_3^2u_5$, $v_2=u_4^2u_5$. Hence, $(x,y,z,t)=(u_3u_4u_5, u_1u_2, u_1u_4^2u_5, u_2u_3^2u_5)$.

By similar analysis, we may solve the remaining two-monomial equations 
$
x^2-y z t = 0
$
and 
$
x_1x_2x_3 - x_4x_5 = 0.
$
The details are left to the reader.

There are three equations left of size $H=12$ that we have not discussed so far. Equation
$$
x+y+x y z = 0
$$
has a family of solutions $(x,y,z)=(0,0,t)$, otherwise $xy \neq 0$ and
$-z=\frac{x+y}{xy}=\frac{1}{x}+\frac{1}{y}$. Hence $|z|\leq 2$. If $z=0$ then $(x,y,z)=(t,-t,0)$. If $|z|=1$ then $(x,y,z)=(2,2,-1)$ or $(-2,-2,1)$. If $|z|=2$ then $(x,y,z)=(1,1,-2)$ or $(-1,-1,2)$.
Equation
$$
x^2-y^2 = zt
$$
can be reduced, after substitution $x-y=x'$ and $x+y=y'$, to the equation $x'y'=zt$, which is equation \eqref{eq:xypzt} with the constraint that either (i) $x'$ and $y'$ are both even, or (ii) $x'$ and $y'$ are both odd. In each case, the solution is a finite union of polynomial families by Proposition \ref{prop:linsub}. Finally, all integer solutions to the equation
$$
x^2+y^2 = zt
$$
can be parametrized as $x=k(uv-wr)$, $y=k(ur+vw)$, $z=k(u^2+w^2)$, $t=k(v^2+r^2)$ for integers $u,v,w,r,k$.
%$x=u_1u_2u_3(u_4u_5+u_6u_7)$, $y=u_1u_2u_3(u_4u_7-u_5u_6)$, $z=u_1^2u_3(u_4^2+u_6^2)$ and $t=u_2^2u_3(u_5^2+u_7^2)$. 
The proof uses the fact that the set of Gaussian integers ${\mathbb Z}[i]=\{a+bi, \, a,b \in {\mathbb Z}\}$, where $i=\sqrt{-1}$, is a unique factorization domain. We then rewrite the equation as $(x+iy)(x-iy)=zt$, and analyse it as a two-monomial equation. See \cite[p. 159]{andreescu2010introduction} for details.

The discussion above can be summarised in the following Proposition.

\begin{proposition}\label{prop:H12}
The set of integer solutions to every polynomial Diophantine equation $P(x_1,\dots,x_n)=0$ of size $H \leq 12$ is a finite union of polynomial families.
\end{proposition}

When moving to equations of size $H=13$, we have two problems. First, there are equations for which it is unclear whether their solution sets are finite unions of polynomial families. We leave them to the reader as the first open question of this paper.

\begin{question}\label{q:h13}
For each of the equations 
$$
yz = x^3+1,
$$
$$
x^2 y = z^2 \pm 1,
$$
$$
zt = x^2 + y^2 \pm 1,
$$
$$
zt = x^2 y+1,
$$
$$
y z t = x^2 \pm 1
$$
and
$$
1 + x_1x_2x_3 + x_4x_5 = 0,
$$
determine whether the set of integer solutions is a finite union of polynomial families.
\end{question} 

The second problem is that there are equations of size $H=13$ whose solution sets are known to be \emph{not} finite unions of polynomial families.
%, see \cite{vaserstein2010polynomial}.
These are the equations
\begin{equation}\label{eq:Pell}
y^2-2x^2=1
\end{equation} 
and
\begin{equation}\label{eq:negPell}
y^2-2x^2=-1,
\end{equation}
which are known as Pell equation and negative Pell equation, respectively. These equations are interesting because if a pair $(x,y)$ of positive integers solves any of them, then ratio $\frac{y}{x}$ is an approximation to $\sqrt{2}$ with a good trade-off between the quality of the approximation and the size of the denominator. For this reason, these equations has been studied since ancient times. It is easy to check that both equations have infinitely many integer solutions. For example, equation \eqref{eq:Pell} has a solution $(x_0,y_0)=(1,0)$, and direct substitution shows that if $(x_n,y_n)$ is a solution, then so is
$$
x_{n+1} = 3 x_n + 2 y_n, \quad y_{n+1} = 4x_n + 3y_n.
$$ 
This gives an infinite sequence of solutions defined by recurrence relations. A bit more difficult to show that this sequence gives \emph{all} solutions with $x\geq 0$ and $y\geq 0$ (and then all integer solutions can be obtained by changing signs), but this is also well-known, see, for example, Theorem 3.2.1 in  \cite{andreescu2015quadratic}. Similarly, all the non-negative solutions to \eqref{eq:negPell} can be obtained starting from $(x_0,y_0)=(1,1)$ and applying the same recurrence relations. As noted in \cite{vaserstein2010polynomial}, these solutions set are too sparse to be described by polynomial families.

This raises the fundamental philosophical issue of what does it mean to ``solve'' a general polynomial Diophantine equation. What if the solution set is infinite, but cannot be described by polynomial families, or recurrence relations, or in any other structured way? To avoid this issue, we may either consider restricted families of equations for which the meaning of ``solve'' is well-defined, or consider general equations but do not aim to find all solutions. The first route is explored in the next section, in which we study equations in $2$ variables.

%This finishes the analysis of all equations of size $H \leq 10$. 
%With $H=11$, we meet some easy equations with $2$ variables that we will discuss in the next section, as well as the equations
%\begin{equation}\label{eq:x2pxp1pyz}
%x^2+x \pm 1-yz=0.
%\end{equation} 
%We leave them to the reader as the first open question of this paper.
%As another example, consider equation
%\begin{equation}\label{eq:x2pxp1pyz}
%x^2+x+1-yz=0
%\end{equation}
%of size $H=11$.  
%
%\begin{question}\label{q:h11}
%Are the solution sets to the equations \eqref{eq:x2pxp1pyz} polynomial families? If not, are they unions of a finite number of polynomial families?
%\end{question} 

\section{Describing all solutions: equations in 2 variables}\label{sec:2varall}

%As we have seen in the previous section, describing the solution set of a $3$-variable equation may be a quite non-trivial problem even with $H\leq 11$. 
In this section, we will restrict our attention to $2$-variable equations, for which many powerful results and techniques are available. We first remark that there is a general algorithm that, given integers $a,b,c,d,e,f$ as an input, solves the general $2$-variable quadratic equation
\begin{equation}\label{eq:quad2var}
ax^2+bxy+cy^2+dx+ey+f=0.
\end{equation} 
The algorithm is implemented online at \cite{Alpern}. For an equation \eqref{eq:quad2var}, it lists all the solutions if there are finitely many of them, and otherwise describes all solutions as a union of polynomial families or in the form of linear recurrence relations.

Based on this, we can exclude all quadratic equations \eqref{eq:quad2var} from further analysis. Together with the previously excluded families, this eliminates all the equations of size $H\leq 13$ with two exceptions:
\begin{equation}\label{eq:ellH13a}
y^2=x^3+1
\end{equation}
and 
\begin{equation}\label{eq:ellH13b}
y^2=x^3-1.
\end{equation}
Equations of the form $y^2=x^3+k$ are known as Mordell's equations and are well studied. It is known that there is a finite number of integer solutions for each $k\neq 0$, and there is an algorithm that, given $k$, outputs all the solutions. See \cite{gebel1998mordell} for the description of the algorithm and for the explicit list of all solutions in the range $|k|\leq 10,000$. 

More generally, there are known practical algorithms for finding integer solutions to the equations in the form
\begin{equation}\label{eq:ellWei}
y^2 + a x y + c y = x^3 + b x ^2 + d x + e
\end{equation} 
under some minor conditions on the integer coefficients $a,b,c,d,e$ that guarantee that the solution set is finite. One such algorithm is implemented in an open-source and free to use computer algebra system SageMath \cite{zimmermann2018computational}, that can be run online at \url{https://sagecell.sagemath.org/}. To solve \eqref{eq:ellWei}, we run the command
\begin{equation}\label{eq:commandellWei}
sage: EllipticCurve([a,b,c,d,e]).integral\_points()
\end{equation}
For example, command
$$
sage: EllipticCurve([0,0,0,0,1]).integral\_points()
$$
returns $[(-1 : 0 : 1), (0 : 1 : 1), (2 : 3 : 1)]$, which means that the only integer solutions to \eqref{eq:ellH13a} are $(x,y)=(-1,0)$, $(0,\pm 1)$, and $(2,\pm 3)$. In a similar way, we find that the only integer solution to \eqref{eq:ellH13b} is $(x,y)=(1,0)$.

This finishes the analysis of the equations of size $H\leq 13$. Starting with $H\geq 14$, we will exclude all the equations in the form \eqref{eq:ellWei}. After this, the only remaining equation of size $H=14$ is
\begin{equation}\label{eq:ellH14}
y^2+yx^2+x=0.
\end{equation} 
This equation is not directly in the form \eqref{eq:ellWei}, but can be easily reduced to it. Indeed, after multiplication by $4y$ and adding $1$ to both sides, we can rewire the equation as $4y^3+(4x^2y^2+4xy+1)=1$, or $(2xy+1)^2=-4y^3+1$. With new variable $z=2xy+1$, this simplifies to $z^2=-4y^3+1$. Now multiply both sides by $16$ to get $(4z)^2=(-8y)^3+16$, or 
\begin{equation}\label{eq:ellH14aux}
Y^2=X^3+16
\end{equation} 
with new variables $Y=4z=4(2xy+1)$ and $X=-8y$. Note that if $x,y$ are integers then so is $X,Y$. Now, command
$$
sage: EllipticCurve([0,0,0,0,16]).integral\_points()
$$
shows that the only integer solutions to \eqref{eq:ellH14aux} are $(X,Y)=(0,\pm 4)$. For these solutions, $y=-X/8=0$ happen to be an integer, and substitution $y=0$ in \eqref{eq:ellH14} returns $x=0$. Hence, $(x,y)=(0,0)$ is the only integer solution to \eqref{eq:ellH14}.

Computer algebra system Maple has a command Weierstrassform that helps to transform a broad range of equations to the form \eqref{eq:ellWei}. In this example, command 
$$
Weierstrassform(x+x^2 y+y^2,x,y,X,Y)
$$
returns $[X^3 - 1/4 + Y^2, y, x*y + 1/2, (-1 + 2*Y)/(2*X), X]$, that shows that the equation can be transformed to $X^3 - 1/4 + Y^2$ after substitutions $X=y$, $Y=xy+1/2$. The remaining steps can be easily done by hand. A combination of Weierstrassform and EllipticCurve commands allows to solve all $2$-variable equations of size $H\leq 15$, and many equations of higher size.

The algorithms we have discussed so far are the special cases of much more general algorithms applicable to much broader classes of $2$-variable equations
\begin{equation}\label{eq:2vargen}
P(x,y)=0.
\end{equation} 
To introduce them, we need a few definitions. A polynomial $P$ with integer coefficients is called absolutely irreducible if it cannot be written as a product $P=P_1\cdot P_2$ of non-constant polynomials, even if we allow complex coefficients. It is known that if $P$ is irreducible over ${\mathbb Q}$ but not absolutely irreducible, then all integer solutions to \eqref{eq:2vargen} can be determined easily, see e.g. \cite{grechuk2021diophantine}. On the other hand, if $P$ is reducible over ${\mathbb Q}$, then equation \eqref{eq:2vargen} reduces to equations of the same form for each of the factors. Hence, we may assume that $P$ in \eqref{eq:2vargen} is absolutely irreducible. In this case, the set of all complex solutions to \eqref{eq:2vargen} forms a connected surface. The \emph{genus} $g$ of such surface is the maximum number of cuttings that can be made along non-intersecting closed simple curves on the surface without making it disconnected. The genus–degree formula 
\begin{equation}\label{eq:gendeg}
g \leq \frac{1}{2}(d-1)(d-2),
\end{equation} 
where $d$ is the degree of $P$, implies that all quadratic polynomials have genus $0$, while all cubic polynomials have genus at most $1$. Poulakis \cite{poulakis1993points, poulakis2002solving} developed practical algorithm to solve all $2$-variable equations of genus $g=0$. The algorithm can decide whether a given equation has finite or infinite number of solutions, list all solutions in the former case, and describes them in the parametric form and/or using recurrence relations in the latter case.

Hence, it suffices to consider equations with $g \geq 1$. In this case, the fundamental theorem of Siegel \cite{siegel1929uber} states that there is always a finite number of integer solutions. The combination of Poulakis and Siegel theorems makes the problem of solving \eqref{eq:2vargen} completely well-defined.

In 1970, Baker \cite{baker1970integer} developed an effective upper bound for the absolute value of all possible solutions to \eqref{eq:2vargen} as an explicit function of the coefficients of $P$, provided that $g=1$. This gives an algorithm to list all the solutions of an arbitrary genus $1$ equation. In particular, by the genus–degree formula \eqref{eq:gendeg}, this result covers all $2$-variable cubic equations. While Baker's bounds are enormous and the corresponding algorithm is impractical, a practical method for finding all integer solutions to genus $1$ equations was later developed by Stroeker and Tzanakis \cite{stroeker2003computing}.

Further, Baker \cite{baker1969bounds} developed in 1969 a general method for solving equations in the form 
\begin{equation}\label{eq:2yPx}
y^2=P(x),
\end{equation} 
where $P(x)$ is a polynomial of arbitrary degree that has at least three simple (possibly complex) zeros. As proved in \cite{grechuk2021diophantine}, this implies the method to determine all integer solutions to the equation 
\begin{equation}\label{eq:quady}
a(x) y^2+b(x)y+c(x)=0.
\end{equation}
where $a(x)$, $b(x)$ and $c(x)$ are arbitrary polynomials with integer coefficients. Indeed, if \eqref{eq:quady} has an integer solution, then $b^2(x)-4a(x)c(x)$ must be a perfect square, and we can apply Baker's algorithm to determine all such $x$, see \cite{grechuk2021diophantine} for details.

This allows us to focus on the equations of degree at least $4$ that are at least cubic in each of the variables. The simplest examples of such equations are, say,  $y^3=x^4+1$ or $y^3=x^4+x+1$. However, such equations are covered by another theorem of Baker \cite{baker1969bounds}, who developed an algorithm for listing all the solutions of the equation 
\begin{equation}\label{eq:y3Px}
y^m=P(x),
\end{equation} 
provided that $m\geq 3$, and $P(x)$ is a polynomial with integer coefficients of degree at least $3$ with at least two simple zeros. In 1984, Brindza \cite{brindza1984s} showed the the conditions on $P$ can be significantly relaxed. It is easy to see \cite{grechuk2021diophantine} that this result also implies the algorithm for solving equations of the form
\begin{equation}\label{eq:ay3Px}
a y^m=P(x).
\end{equation} 
Baker's and Brindza's methods for solving equations \eqref{eq:2yPx} and \eqref{eq:y3Px} are impractical even for the equations with small coefficients. However, there are practical methods for which we do not have proof that they work in general, but which seem to work for any individual equation in this form. For example, Bruin and Stoll \cite{bruin2008deciding} decided the solvability in rationals of all the equations \eqref{eq:2yPx} where $P$ is a square-free, has degree at most $6$, and has integral coefficients of absolute value at most $3$. More recently, Hashimoto and Morrison \cite{hashimoto2020chabauty} determined the set of all rational solutions for a large family of the equations in the form \eqref{eq:y3Px}.

Based on this, we eliminate equations of the form \eqref{eq:ay3Px} from further analysis. After this, the smallest non-eliminated ones are 
\begin{equation}\label{eq:2varH26}
x^3y+y^3\pm x = 0
\end{equation}
of size $H=26$ and the next-smallest are
\begin{equation}\label{eq:2varH27}
x^3y+y^3\pm x + 1 = 0
\end{equation}
of size $H=27$. These equations can be easily solved directly (exercise!), but we instead note that all such equations are covered by another deep result. Let ${\cal F}$ be a family of polynomials
$$
P(x,y) = \sum_{i=0}^m \sum_{j=0}^n a_{ij} x^i y^j
$$
with integer coefficients $a_{ij}$ of degree $m>0$ in $x$ and $n>0$ in $y$ which are irreducible over ${\mathbb Q}[x,y]$, and such that either
\begin{itemize}
\item[(C1)] there exists a coefficient $a_{ij}\neq 0$ of $P$ such that $ni+mj>mn$,

or 

\item[(C2)] the sum of all monomials $a_{ij}x^iy^j$ of $P$ for which $ni + mj = nm$ can be decomposed into a product of two non-constant relatively prime polynomials in ${\mathbb Z}[x,y]$.
\end{itemize}
In 1887, Runge \cite{runge1887ueber} proved if $P\in {\cal F}$ then equation $P=0$ has at most finitely many integer solutions. In 1992, Walsh \cite{walsh1992quantitative} developed an effective upper bound for the size of possible solutions, which implies the existence of an algorithm for listing all the solutions. See \cite{beukers2005implementation} for a practical implementation of this algorithm. Note that equations \eqref{eq:2varH26} and \eqref{eq:2varH27} satisfy (C1), because in this case $n=m=3$, there is a non-zero coefficient $a_{31}=1$, and, for $i=3$ and $j=1$, we have $ni+mj=3\cdot 3 + 3\cdot 1 > 3 \cdot 3 = mn$. Walsh's theorem allows us to exclude all equations that satisfy either (C1) or (C2) from further analysis.

This excludes all equations of size $H\leq 27$, and all equations of size $H=28$ with three exceptions:
\begin{equation}\label{eq:2varH28a}
x^4+x y+y^3=0,
\end{equation}
\begin{equation}\label{eq:2varH28b}
y^3+y=x^4+x
\end{equation}
and
\begin{equation}\label{eq:2varH28c}
y^3-y=x^4-x.
\end{equation}
%\begin{equation}\label{eq:2varH29a}
%x^4+x y+y^3+1=0
%\end{equation}
%and
%\begin{equation}\label{eq:2varH29b}
%x^4+x y+y^3-1=0.
%\end{equation}
The listed equations does not satisfy (C1) because for them we have $m=4$, $n=3$, and there is no non-zero coefficient $a_{ij}$ with $3i+4j>12$. Equality $3i+4j=12$ holds for coefficients $a_{40}$ and $a_{03}$, and polynomials $x^4\pm y^3$ are irreducible, hence (C2) also fails. 

Equation \eqref{eq:2varH28a} has genus $2$. Computer algebra system Magma has built-in  method (called Chabauty) for finding all rational solutions for some genus $2$ equations, and the method happens to work for this particular equation, returning that the only rational solution is $x=y=0$. Equations \eqref{eq:2varH28b} and \eqref{eq:2varH28c} have genus $3$, and this Magma function is not applicable to them. By Siegel1's Theorem \cite{siegel1929uber}, the sets of their integer solutions are finite. The direct search returns solutions $(x,y)=(-1,0)$, $(0,0)$, and $(1,1)$ for \eqref{eq:2varH28b} and $(x,y)=(0,-1)$, $(0,0)$, $(0,1)$, $(1,-1)$, $(1,0)$ and $(1,1)$ for \eqref{eq:2varH28c}, but the problem is to prove that no other solutions exists. We leave this to the reader as open questions.

\begin{question}\label{q:h28a}
Find all integer solutions to \eqref{eq:2varH28b}.
\end{question}

\begin{question}\label{q:h28b}
Find all integer solutions to \eqref{eq:2varH28c}.
\end{question}

\section{Finding the solution set if it is finite}\label{sec:fin}

%\subsection{Some trivial methods}

Now let us return to Diophantine equations in $3$ or more variables. In this case the problem of ``solving'' the equation is, in general, not well posed: if the solution set is infinite but not a finite union of polynomial families and cannot be described by recurrence relations, then was counts as an ``acceptable description'' of this solution set? For the sets with no obvious ``structure'' this problem is more philosophical than mathematical, and we will not discuss it further. Instead, we will focus on the following problem, which is completely well-defined. 

\begin{problem}\label{prob:fin}
Given a polynomial Diophantine equation, decide whether its solution set is finite, and if so, list all the solutions.
\end{problem}

Note that proving that the solution set is infinite completely solves Problem \ref{prob:fin}, and no further analysis is required. In Section \ref{sec:polfam}, we have completed the analysis of all equations of size $H \leq 12$, so we may move to $H=13$.
% If we focus on Problem \ref{prob:fin} only, equations \eqref{eq:x2pxp1pyz} in the Open Question \ref{q:h11} are trivial, because they obviously has infinitely many integer solutions. Indeed, we may choose $z=1$, take arbitrary $x$, and set $y=x^2+x\pm 1$. This finishes the analysis of all equations of size $H \leq 11$.
We first remark that there are many equations like
\begin{equation}\label{eq:almostPyth}
x^2+y^2=z^2 + 1
\end{equation}
that has infinitely many integer solutions for some fixed value of one of the variables (for \eqref{eq:almostPyth}, take $y=1$ and $x=z$). For such equations, Problem \ref{prob:fin} is trivial hence they may be excluded. After this, the only non-excluded equation of size $H\leq 13$ is
\begin{equation}\label{eq:nearlyPyth}
x^2+y^2=z^2 - 1
\end{equation}
This equation has at most a finite number of integer solutions for any fixed $x$, $y$, or $z$, but still has infinitely many integer solutions. Indeed, for any integer $t$, we have a solution $x=2t^2$, $y=2t$, $z=2t^2+1$. Integer solutions to \eqref{eq:almostPyth} and \eqref{eq:nearlyPyth} are known as ``almost polynomial triples'' and ``nearly polynomial triples'', respectively. See \cite{frink1987almost} for the complete description of the solution sets to these equations.
% see Frink, O.: Almost Pythagorean Triples, Mathematics Magazine, Vol. 60, No. 4, October 1987, pp. 234-236 describes them all.

More generally, for any integer $a$, equation
\begin{equation}\label{eq:aPyth}
x^2+y^2=z^2 + a
\end{equation}
has infinitely many integer solutions. Indeed, we can rewrite the equation as $x^2-a=y^2-z^2=(y-z)(y+z)$. For simplicity, assume that $y-z=1$, so that $x^2-a=y+z=2z+1$, from which we can find $z=(x^2-a-1)/2$. Now, if $a=2k-1$ is odd, take $x=2t$, $z=2t^2-k$, and $y=z+1=2t^2-k+1$, while if $a=2k$ is even, take $x=2t+1$, $z=2t^2+2t-k$, and $y=z+1=2t^2+2t-k+1$. As a side note, we remark that a complete description of the solution set to \eqref{eq:aPyth} with $a=3$ has been an open question until Vaserstein \cite[Example 15]{vaserstein2010polynomial} proved that it is the union of two polynomial families.

%Next, equation
%\begin{equation}\label{eq:PythH14a}
%x^2 + y^2 = z^2 + 2
%\end{equation}
%has at most a finite number of integer solutions for any fixed $x$, $y$, or $z$, but has infinitely many integer solutions in the form $x=2t+1$, $y=2t^2+2t-1$, $z=2t^2+2t$. Similarly, equation 
%\begin{equation}\label{eq:PythH14b}
%x^2 + y^2 = z^2 - 2
%\end{equation}
%$x=2t+1$, $y=2t^2+2t-1$, $z=2t^2+2t$. As another example, equation 
%\begin{equation}\label{eq:PythH15}
%x^2 + y^2 = z^2 + 3
%\end{equation}
%has infinitely many integer solutions because we may take $x=2t$, $y=2t^2-2$, and  $z=2t^2-1$. As a side note, we remark that a complete description of the solution set to \eqref{eq:PythH15} has been an open question until Vaserstein \cite[Example 15]{vaserstein2010polynomial} proved that it is the union of two polynomial families.

More generally, if there exist polynomials $Q_1(t), \dots, Q_n(t)$ with integer coefficients, not all constant, such that
\begin{equation}\label{eq:PQ1On}
P(Q_1(t), \dots, Q_n(t)) \equiv 0
\end{equation}
then equation $P(x_1,\dots, x_n)=0$ has infinitely many integer solutions. In general, deciding the existence of such polynomials is a quite non-trivial problem. However, we can at least verify \eqref{eq:PQ1On} for polynomials with small degree and coefficients, and exclude the equations for which we managed to find the corresponding $Q_i$. This method allowed us to solve all the remaining equations of size $H \leq 16$.

%\subsection{Sun of two squares}

The first equation of size $H=17$ we discuss is 
\begin{equation}\label{eq:H17a}
y^2 + z^2 = x^3 - 1
\end{equation}
We will prove that it has infinitely many integer solutions. Equivalently, there exists infinitely many integers $x$ such that $x^3-1$ is the sum of two squares. Identity
$$
(a^2+b^2)(c^2+d^2)=(ac+bd)^2+(ad-bc)^2
$$
shows that if two integers can be represented as sum of squares, then so is their product. Because $x^3-1=(x-1)(x^2+x+1)$, it suffices to find $x$ such that $x+1=u^2$ and $x-1=v^2+w^2$, or $v^2+w^2=u^2-2$. But this is equation \eqref{eq:aPyth} with $a=-2$ which has infinitely many integer solutions. In fact, this leads to explicit parametric family $x=3 + 8 t + 12 t^2 + 8 t^3 + 4 t^4$, $y=5 + 20 t + 38 t^2 + 40 t^3 + 24 t^4 + 8 t^5$ and $z=-1 - 8 t - 28 t^2 - 44 t^3 - 44 t^4 - 24 t^5 - 8 t^6$ of the solutions to \eqref{eq:H17a}, but the coefficients in these polynomials are too large for the direct search. This is the reason why \eqref{eq:H17a} has not been excluded automatically and required explicit argument.
%While verification is straightforward, one may ask how to come up with these formulas. The idea is that $x^3-1=(x-1)(x^2+x+1)$. If two integers can be represented as sum of squares, then so is their product. So, it suffices if $x+1=u^2$ and $x-1=a^2+b^2$, or $a^2+b^2=u^2-2$. And this equation (of size $H=14$) is easier and its parametric solution can be found by direct search.  

Another equation of size $H=17$ that requires attention is
\begin{equation}\label{eq:H17b}
y(x^2-y) = z^2+1
\end{equation}
We will prove that it has no integer solutions. For this, we will need a well-known fact \cite{grechuk2021diophantine} that all odd prime factors of a sum of squares $z^2+1$ must be congruent to $1$ modulo $4$. Hence the same is true for the odd prime factors of positive integers $y$ and $x^2-y$. Because the product of any number of such primes is again $1$ modulo $4$, this implies that if $z^2+1$ is odd, then both $y$ and $x^2-y$ are congruent to $1$ modulo $4$, but then $x^2$ is congruent to $2$ modulo $4$, a contradiction. If $z^2+1$ is even, its prime factorization contains exactly one factor of $2$, which goes to either $y$ or $x^2-y$, resulting in $x^2$ being $3$ modulo $4$, again a contradiction.

This finishes the analysis of equations of size $H=17$. For $H=18$, we start with the equation similar to \eqref{eq:H17a},
\begin{equation}\label{eq:H18a}
y^2+z^2 = x^3 - 2
\end{equation}
Unlike $x^3-1$, $x^3-2$ does not factorise, so different technique is required. We will present a solution given by Max Alekseyev in the comment to mathoverflow question\footnote{\url{https://mathoverflow.net/questions/409857/representing-x3-2-as-a-sum-of-two-squares}}.
Let $x=t^2+2$ for some integer $t$. Then $x^3-2=(t^2+2)^3-2=(t^3+3t)^2+(3t^2+6)$. It is left to select $y=t^3+3t$ and note that equation $z^2 = 3t^2 + 6$ has infinitely many integer solutions. This can be checked directly or using the Gauss theorem \cite[p. 57]{mordell1969diophantine}, that states that the general quadratic equation with integer coefficients
$$
az^2 + bzt + ct^2 + dz + et + f = 0
$$
such that $D=b^2-4ac > 0$, $D$ is not a perfect square, and $\Delta = 4acf + bde - ae^2 - cd^2 - fb^2 \neq 0$ has either no integer solutions or infinitely many of them. In our case, $D=-4\cdot 3 \cdot (-1) = 12>0$, $\Delta \neq 0$, and there is an integer solution, say $z=3$, $t=1$. This finishes the proof.

Another equation of size $H=18$ is
\begin{equation}\label{eq:H18b}
x^2+y^2+x y z -2 = 0
\end{equation}
This equation requires completely different techniques called Vieta jumping. The idea is that if $(x,y,z)$ is any solution to \eqref{eq:H18b}, then $t=x$ is a solution to the quadratic equation 
$$
t^2+y z t +y^2 -2 = 0,
$$
and this equation has another solution $t'=-y z - x = \frac{y^2-2}{x}$, which is also an integer. Hence, any solution $(x,y,z)$ to \eqref{eq:H18b} produces another solution $(-yz-x,y,z)$, and, by a similar argument, one more solution $(x,-xz-y,z)$. The technique suggests to consider a solution with $|x|+|y|+|z|$ minimal and either prove that there is no such solution, or find such minimal solution and then produce infinitely many other solutions by the transformations above.

Returning to equation \eqref{eq:H18b}, it has solutions $(x,y,z)=(\pm 1, \pm 1, 0)$ with $z=0$. Let us prove that there are no other solutions. Assume the contrary and let $(x,y,z)$ be a solution with $z\neq 0$ such that $|x|+|y|+|z|$ is minimal. By symmetry, we may assume that $|x|\geq |y|$. By Vieta jumping, there is another solution $(x',y,z)$ with $|x'|\geq |x|$ and $xx' = y^2-2$. Then $|y^2-2| = |x||x'| \geq x^2 \geq y^2$, hence either $y^2=0$ or $y^2=1$. But $y^2=0$ implies $x^2=2$ which is impossible, while $y^2=1$ implies that $z=0$, a contradiction. Hence, the only solutions are $(x,y,z)=(\pm 1, \pm 1, 0)$.

To apply this technique to general equation $P(x_1,\dots,x_n)=0$, let us denote $S$ the set of all variables $x_i$ for which the equation can be written as
$$
a_i x_i^2 + Q_i x_i + R_i = 0
$$
where $|a_i|=1$ and $Q_i$ and $R_i$ are polynomials in other variables. We can then use any computer algebra system to solve optimization problem of maximizing $t$ over $(x_1,\dots,x_n,t)\in {\mathbb R}^{n+1}$ subject to constraints $P=0$, $|x_i|\geq t$ for each $i$, and $|-(Q_i/a_i)-x_i|\geq |x_i|$ for each variable $x_i \in S$. If the optimal value $t^*$ to this optimization problem is infinite, then the method does not work for this equation. But if $t^*<\infty$, then we have $\min\{|x|,|y|,|z|\}\leq t^* < \infty$ for any solution $(x,y,z)$ with $|x|+|y|+|z|$ minimal. So, we next check, for each integer $t$ such that $|t|\leq t^*$ and each $i=1,\dots,n$, whether the equation has any integer solutions with $x_i=t$. If there are no such solutions, then the equation has no integer solutions at all. If there are such solutions, we next check whether any of them produces an infinite chain of solutions via Vieta jumping. In the rest of the paper, we will exclude the equations solvable by this method.

This finishes the analysis of all equations of size $H \leq 18$. Starting with $H=19$, we get several equations in the form
$$
ax^2+bx+c+dy+exyz=0,
$$ 
such as, for example, 
$
1+x-x^2+2 y+x y z = 0,
$
$
1+x^2+3 y+x y z = 0,
$
etc. These equations are not covered by techniques discussed so far, but are in fact easy. Indeed, solutions with $x=0$ or $y=0$ can be found by direct substitution. Otherwise fraction $\frac{dy+c}{x}$ is an integer, hence we have $|dy+e|\geq |x|$, or $|x/y|\leq |d+e/y|\leq |d|+|e|$. Then $|z|=|(ax^2+bx+c+dy)/(exy)|\leq|(ax^2+bx+c+dy)/(xy)|=|ax/y+b/y+c/(xy)+d/x| \leq |a|(|d|+|e|)+|b|+|c|+|d|$. Now direct search in this region returns the full set of solutions.

A more interesting equation that require a new idea is
\begin{equation}\label{eq:H18c}
y^2+z^2=x^3+3
\end{equation}
To solve it, recall that a positive integer is the sum of two squares if and only if all its prime factors congruent to $3$ modulo $4$ enters its prime factorization an even number of times, see e.g. \cite{grechuk2021diophantine}. In particular, this implies that if $a$ and $b$ do not share prime factors congruent to $3$ modulo $4$, and $ab$ is the sum of two squares, then so are both $a$ and $b$. Now note that $x(x^3+3)=(x^2-1)^2+(2x^2+3x-1)$. Let $x$ be any integer such that $2x^2+3x-1$ is a perfect square (there are infinitely many such integers). Then $x$ is not divisible by $3$ (otherwise $2x^2+3x-1$ would be $2$ mod $3$ and could not be a perfect square), hence $x$ and $x^3+3$ are co-prime. But their product $x(x^3+3)$ is the sum of two squares, hence $x^3+3$ is a sum of two squares as well. 

The same method allows to solve many other similar equations, such as 
$
y^2+z^2 = x^3+x+1,
$
$
y^2+z^2 = x^3-x-1,
$
$
y+y^2+z^2 = x^3-1,
$
etc. (The last equation after multiplication by $4$ can be rewritten as $(2y+1)^2+(2z)^2=4x^3-3$, so it suffices to prove that $4x^3-3$ is the sum of squares infinitely often, and then the same method applies). We will exclude any further equations solvable by this method. This finishes the analysis for $H \leq 19$.

The only new equations of size $H=20$ are homogeneous quadratic equations like
$$
x^2+y^2=3z^2.
$$
The only integer solution to this equation is $(x,y,z)=(0,0,0)$. Indeed, if there is any other solution then we can divide it by any common factor and obtain a new solution for which $(x,y,z)$ are co-prime. However, the sum of squares $x^2+y^2$ is divisible by $3$ only if both $x$ and $y$ are divisible by $3$. But in this case $x^2+y^2$ is divisible by $9$, hence $z$ is divisible by $3$, a contradiction with the co-primality assumption. Famous Hasse–Minkowski theorem (Hasse principle) states that if a homogeneous quadratic equation has non-zero real solutions but no non-zero integer solutions, then this can always be proved by divisibility analysis modulo some $p$ as above. This allows us to exclude such equations as well. 

For $H=21$, the only equation of different type is
\begin{equation}\label{eq:h21}
y(x^2+2) = 2z^2-1.
\end{equation}
So far we have used only the information about prime factors of sum of two squares, while this equation requires the analysis of prime factors of other quadratic polynomials, in this case $x^2+2$ and $2z^2-1$. As shown in \cite{grechuk2021diophantine}, all odd prime factors of $x^2+2$ must be $1$ or $3$ modulo $8$, while all prime factors of $2z^2-1$ must be $1$ or $7$ modulo $8$. A combination of these facts imply that if $(x,y,z)$ solves \eqref{eq:h21}, then all prime factors of $x^2+2$ are congruent to $1$ modulo $8$. But then $x^2+2$ must be itself congruent to $1$ modulo $8$, which is a contradiction. We refer to \cite{grechuk2021diophantine} how to apply this method in general, but here will not list further equation solvable in this way.

For $H=22$, we start to meet equations like
\begin{equation}\label{eq:h21a}
y^2+y z+z^2 = x^3-x
\end{equation}
that require the analysis of which integers can be represented in the form $y^2+yz+z^2$. Let $S$ be the set of all such integers. It is known that $S$ is also the set of integers representable as $3y^2+z^2$, and also the set of integers $n$ such that every prime $p$ of the form $p=3k+2$ enters the prime factorization of $n$ in the even power. We need to prove that $x^3-x$ belongs to $S$ for infinitely many $x$. Choose any odd $x$ such that $2x^2 - 2x-4 = 3t^2$ for some integer $t$ (there are infinitely many such $x$). Then $(x^3-x)(x-2)=(x^2-x-2)^2+(2x^2 - 2x-4)$ belongs to $S$. Because $(x^3-x)$ and $(x-2)$ do not share any prime factors in the form $p=3k+2$, this implies that $x^3-x \in S$. The same method allows to solve other equations of this type, such as 
$
y^2+y z+z^2=x^3+x
$
and
$
y^2+y z+z^2=x^3-2.
$

The next equation we discuss is
$$
y(z^2-y) = x^3+2. 
$$
We present a solution given by Mathoverflow user Tomita\footnote{\url{https://mathoverflow.net/questions/411958}}. By considering the equation as quadratic in $y$, we conclude that it has infinitely many integer solutions if and only if the determinant $D=(-z^2)^2-4(x^3+2)=z^4-4x^3-8$ is a perfect square infinitely often. Now assume that $x=-3t^2-2t-2$ and $z=3t+1$ for some integer $t$. Then $D=(3t+1)^4-4(-3t^2-2t-2)^3-8=(12t^2+8t+25)(3t^2+2t+1)^2$. It is left to remark that $12t^2+8t+25$ is a perfect square for infinitely many integers $t$. 

The same method solves another equation, 
$$
xyz=x^3+y^2-2.
$$ 
We need $D=x^2z^2-4x^3+8$ to be a perfect square. Select $x=6t^2+1$ and $z=6t$, then $D=4(6t^2-1)^2(3t^2+1)$. It is left to note that there are infinitely many integers $t$ such that $3t^2+1$ is a perfect square.

However, we currently do not see how to use these (or other) methods to solve similar equations
\begin{equation}\label{eq:h21b}
y(z^2-y) = x^3-2
\end{equation}
and
\begin{equation}\label{eq:h21c}
xyz=x^3+y^2+2.
\end{equation}
These are the only remaining open equations of size $H \leq 22$. A computer search for polynomials $x=Q(t)$ and $z=R(t)$ with small degree and coefficients returns no polynomials for which the same method works. Hence, we need either a deeper search for polynomials with large coefficients, or a new idea. We will leave these equations to the readers as open questions.

\begin{question}\label{q:h22a}
Are there infinitely many integer solutions to \eqref{eq:h21b}?
\end{question}

\begin{question}\label{q:h22b}
Are there infinitely many integer solutions to \eqref{eq:h21c}?
\end{question}

One may also study Problem \ref{prob:fin} for some restricted families of polynomials. For example, if we restrict the number of variables and consider $2$-variable equations only, then the smallest equations for which Problem \ref{prob:fin} is open are equations \eqref{eq:2varH28b} and \eqref{eq:2varH28c} of size $H=28$, see open questions \ref{q:h28a} and \ref{q:h28b}. 

Another nice class of equations we may consider are symmetric equations, ones that are invariant after cyclic shirt of variables. The smallest symmetric equation not directly solvable by the methods described above turns out to be the equation 
\begin{equation}\label{eq:sym}
x^2y+y^2z+z^2x=1
\end{equation}
of size $H=25$.

\begin{question}\label{q:sym}
Solve Problem \ref{prob:fin} for the equation \eqref{eq:sym}.
\end{question}

Finally, we may also restrict the number of monomials. It is easy to solve all $2$-monomial equations \cite{grechuk2021diophantine}, hence the first interesting case is $3$-monomial ones. The smallest $3$-monomial equation which seems to be not solvable by the described methods is the equation
\begin{equation}\label{eq:3mon}
x^3y^2=z^3+2
\end{equation}
of size $H=42$. This equation has obvious solutions $(x,y,z)=(1,\pm 1, -1)$. Note that any integer $n$ can be represented in the form $x^3y^2$ if and only if for every prime number $p$ dividing $n$, $p^2$ also divides $n$. Such integers are called powerful numbers. So, the question is to find all integers $z$ such that $z^3+2$ is a powerful number.

\begin{question}\label{q:3mon}
Find all integer solutions to the equation \eqref{eq:3mon}.
\end{question}

\section{Existence of solutions: Hilbert 10th problem}\label{sec:Hilbert}

In addition to Problem \ref{prob:fin}, one may consider the following problem with Yes/No answer.

\begin{problem}\label{prob:main}
Given a polynomial Diophantine equation, check whether it has any integer solution.
\end{problem}

Hilbert's 10th problem asks for a general method for solving Problem \ref{prob:main} for all Diophantine equations. Building on the work of Davis, Putnam and Robinson \cite{davis1961decision}, Matiyasevich \cite{matijasevic1970enumerable} proved in 1970 that no such general algorithm exists. 
%See excellent recent surveys of Gasarch \cite{gasarch2021hilbertb, gasarch2021hilberta, gasarch2021hilberts} for a detailed discussion for which families of Diophantine equations Problem \ref{prob:main} can be solved, and for which families it is known to be undecidable. 
%See excellent recent surveys of Gasarch \cite{gasarch2021hilbertb, gasarch2021hilberta, gasarch2021hilberts} for a detailed discussion for which families of Diophantine equations can be solved, and for which families it is known to be undecidable.
See excellent recent surveys of Gasarch \cite{gasarch2021hilbertb, gasarch2021hilberta, gasarch2021hilberts} for a detailed discussion for which Diophantine equations Problem \ref{prob:main} can be solved, and in which cases it is known to be undecidable.
For all the equations we left open in the previous sections, Problem \ref{prob:main} is trivial because these equations have some obvious small solutions.

In \cite{grechuk2021diophantine}, we found the smallest Diophantine equation for which Problem \ref{prob:main} is currently open. This is the equation
\begin{equation}\label{eq:h31main}
y(x^3-y)=z^3+3
\end{equation}
of size $H=31$.

\begin{question}\label{qu:smallest}
Do there exist integers $x,y,z$ satisfying \eqref{eq:h31main}? 
\end{question}

The same question can also be asked for restricted families of equations. Among the $2$-variable equations, the smallest open are
\begin{equation}\label{eq:H32a}
y^3+xy+x^4+4=0,
\end{equation}
\begin{equation}\label{eq:H32b}
y^3+xy+x^4+x+2=0,
\end{equation}
\begin{equation}\label{eq:H32c}
y^3+y=x^4+x+4
\end{equation}
and
\begin{equation}\label{eq:H32d}
y^3-y=x^4-2x-2
\end{equation}
of size $H=32$.

\begin{question}\label{qu:2var}
Determine whether each of the equations \eqref{eq:H32a}-\eqref{eq:H32d} have any integer solution.
\end{question}

%The smallest open symmetric equation is
%\begin{equation}\label{eq:H37sym3}
%x^3 + y^3 + z^3 + xyz = 5
%\end{equation}
%of size $H=37$. 
%
%\begin{question}
%Do there exist integers $x,y,z$ satisfying \eqref{eq:H37sym3}? 
%\end{question}

The smallest open symmetric equation is
\begin{equation}\label{eq:H39sym3}
x^3+x+y^3+y+z^3+z = x y z + 1
\end{equation}
of size $H=39$. 

\begin{question}
Do there exist integers $x,y,z$ satisfying \eqref{eq:H39sym3}? 
\end{question}

Finally, the smallest open $3$-monomial equation is
\begin{equation}\label{eq:H46term3}
x^3y^2 = z^3 + 6
\end{equation}
of size $H=46$.
 
\begin{question}
Do there exist integers $x,y,z$ satisfying \eqref{eq:H46term3}? 
\end{question}

In addition, we may consider the smallest open equations with respect to alternative measures of size. As noted in the introduction, a natural measure of ``length'' of a monomial $M$ of degree $d$ with coefficient $a$ is $l(M)=\log_2 |a|+d$. Then we can define the length $l(P)$ of a polynomial $P$ consisting of $k$ monomials with coefficients $a_1, \dots, a_k$ and degrees $d_1, \dots, d_k$, respectively, as
\begin{equation}\label{eq:ldef}
l(P)=\sum_{i=1}^k (\log_2 |a_i|+d_i).
\end{equation}
Then, instead of ordering the equations by $H$, we may order them by length $l$, or, equivalently, by an integer
$$
L(P) := 2^{l(P)} = \prod_{i=1}^k|a_i|\cdot 2^{\sum_{i=1}^k d_i}.
$$
Note that the formula for $L(P)$ is the same as the formula \eqref{eq:Hdef} for $H(P)$, except that the summation is replaced by a product.

As established in \cite{grechuk2021diophantine}, the shortest equations for which Problem \ref{prob:main} is open are the equations
\begin{equation}\label{eq:l10main1}
y(x^3-y) = z^4+1,
\end{equation}
\begin{equation}\label{eq:l10main2}
2 y^3 + x y + x^4 + 1 = 0
\end{equation}
and 
\begin{equation}\label{eq:l10main3}
x^3 y^2 = z^4+2
\end{equation}
that have length $l=10$. 

\begin{question}\label{qu:shortest}
Do there exist integers $x,y,z$ satisfying \eqref{eq:l10main1}? 
\end{question}

\begin{question}\label{qu:shortest2}
Do there exist integers $x,y$ satisfying \eqref{eq:l10main2}? 
\end{question}

\begin{question}\label{qu:shortest3}
Do there exist integers $x,y,z$ satisfying \eqref{eq:l10main3}? 
\end{question}

\section{Conclusions}\label{sec:concl}

We have ordered all Polynomial Diophantine equations by a parameter $H$ defined in \eqref{eq:Hdef} and tried to solve the equations in that order. We have considered the following problems in the decreasing level of difficulty.

\begin{itemize}

\item Completely solve the equation: list all solutions if there are finitely many and describe all solutions (for example, as a union of polynomial families) if the solution set is infinite.

\item Determine whether the solution set is finite, and if yes, list all the solutions.

\item Check whether an equation has any integer solution.

\end{itemize}

For each of the problems, we have identified the smallest equations for which the problem is open. In some cases, we also identified the smallest open equations is certain families, such as the smallest open $2$-variable, symmetric, or $3$-monomial equations. 
The list of current smallest open equations can also be found on Mathoverflow \cite{BG2021,BG2021b,BG2021c}, where the plan is to always keep the list up-to-dated.

%In the Mathoverflow questions 
%\footnote{\url{https://mathoverflow.net/questions/400714/can-you-solve-the-listed-smallest-open-diophantine-equations}} 
%\footnote{\url{https://mathoverflow.net/questions/411958/on-the-smallest-open-diophantine-equations-beyond-hilberts-10-problem}}

%\begin{center}
%\begin{tabular}{ |c|c|c| } 
% \hline
% Equation & Size & Comment \\ 
% \hline\hline
% $y(x^3-y)=z^3+3$ & $H=31$ & The smallest (in $H$) open equation \\ 
%\hline
% $y^3+xy+x^4+4=0$ & $H=32$ & The smallest open $2$-variable equations \\ 
% $y^3+xy+x^4+x+2=0$ & $H=32$ &  \\ 
% $y^3+y=x^4+x+4$ & $H=32$ & \\ 
% $y^3-y=x^4-2x-2$ & $H=32$ & \\ 
%\hline
% $x^3 + y^3 + z^3 + xyz = 5$ & $H=37$ & The smallest open symmetric equation \\ 
%\hline
% $x^3y^2 = z^3 + 6$ & $H=46$ & The smallest open $3$-monomial equation \\ 
%\hline
% $y(x^3-y) = z^4+1$ & $l=10$ & The shortest (in $l$) open equations \\ 
% $2 y^3 + x y + x^4 + 1 = 0$ & $l=10$ &  \\ 
% $x^3 y^2 = z^4+2$ & $l=10$ & \\ 
% \hline
%\end{tabular}
%\end{center}

We suspect that some of the open equations listed in this paper are relatively easy, and are suitable for the first research project of a graduate or even undergraduate student. On the other hand, we are confident that some of our equations are quite difficult and may stimulate the development of new methods and techniques in number theory. 

%Such equations may serve as nice test cases for any known or new technique for determining integer and/or rational points on curves and surfaces. 
%
%To taste how difficult Problem \ref{prob:allsol} can be, think how you would describe \emph{all} integer solutions to the simple equation $x_1x_4-x_2x_3=1$ of size $H=9$. In 2010, Vaserstein \cite{vaserstein2010polynomial} found a way how to do this using polynomial expressions with $46$ parameters. See Section 10.4 of \cite{grechuk2019theorems} for an accessible description of this theorem.

%You may also try to solve any of Problems \ref{prob:main}-\ref{prob:allsol} for any restricted family of polynomials, again ordered by $H$. For example, you may consider:
%\begin{itemize}
%\item Homogeneous polynomials. For this family, Problem \ref{prob:main} is trivial because of zero solution, but Problems \ref{prob:nonzero}-\ref{prob:allsol} all make sense.
%\item Polynomial in specific number of variables, for example, in $2$ variables.
%\item Polynomials of a given degree, for example, cubic ones.
%\item Symmetric polynomials.
%\item Various intersections of the above families. For example, can you find the simplest (as measured by $H$) symmetric polynomial in $2$ variables for which Problem \ref{prob:finlist} is an open question?
%\end{itemize}

\section*{Acknowledgments}

I thank anonymous mathoverflow user Zidane whose question \cite{Z2018} inspired this work. I also thank other mathoverflow users, including but not limited to  Will Sawin, Fedor Petrov, Andrew R. Booker, Victor Ostrik, Jeremy Rouse, Max Alekseyev and Tomita for very helpful discussions on this topic and for solving some equations in this project that was previously listed as open. I thank Aubrey de Grey for writing his own version of the computer program for enumerating equations, which provides an independent validation that none of the equations of small size has been accidentally missing. I also thank William Gasarch for the interest to publish this paper in SIGAST NEWS Open Problem Column and for proofreading an earlier draft of this paper and providing me with helpful list of comments and suggestions.

\end{document}